\documentclass[12pt]{article}
\usepackage{amssymb,amsmath,amsthm,amsfonts}
\usepackage{url}
\usepackage{multirow}
\usepackage{graphicx}
\usepackage{fancyhdr}
\newtheorem{Theorem}{Theorem}[section]

\newtheorem{assumption}{Assumption}[section]

\newenvironment{tightenumerate}{\begin{list}{--~}{
  \topsep=0.3ex \itemsep=0.3ex \labelsep=1em \parsep=0em
  \listparindent=0em \itemindent=0em
  \settowidth{\labelwidth}{--~} \leftmargin=3.5em
}}{\end{list}}

\newcommand{\imp}{{\tt imp}}
\newcommand{\R}{\mathbb{R}}
\newcommand{\Rtol}{r_{\tt tol}}
\newcommand{\Gtol}{\varepsilon_{\tt tol}}
\newcommand{\Dtol}{\Delta_{\tt tol}}

\newcommand{\DFPP}{{\tt DFPP}}
\newcommand{\lin}{{\tt S_n}}
\newcommand{\bilin}{{\tt S_{2n}}}
\newcommand{\mean}{{\tt S_{n^2/4}}}
\newcommand{\qua}{{\tt S_{n^2/2}}}

\title{\textbf{Adaptive Interpolation Strategies in Derivative-Free Optimization: a case study}}

\author{W.~Hare\thanks{Mathematics, Irving K. Barber School, University of British Columbia, Okanagan Campus, Kelowna, B.C., V1V 1V7, Canada: {\tt warren.hare@ubc.ca}.} \and M.~Jaberipour\thanks{Mathematics and Computer Science, Amirkabir University of Technology, Tehran, Iran: {\tt majid.jaberipour@gmail.com}.}}

\begin{document}
\renewcommand{\topfraction}{.9}
\renewcommand{\textfraction}{0.01}

\maketitle

\begin{abstract}
Derivative-Free optimization (DFO) focuses on designing methods to solve optimization problems without the analytical knowledge of gradients of the objective function. There are two main families of DFO methods: model-based methods and direct search methods.  In model-based DFO methods, a model of the objective function is constructed using only objective function values, and the model is used to guide the computation of the next iterate.  Natural questions in this class of algorithms include how many function evaluations should be used to construct the model?  And, should this number be fixed, or adaptively selected by the algorithm? In this paper, we numerically examine these questions, using Hare and Lucet's Derivative-Free Proximal Point (\DFPP) algorithm \cite{W.L.Hare2014} as a case study.  Results suggest that the number of function evaluations used to construct the model has a huge impact on algorithm performance, and adaptive strategies can both improve and hinder algorithm performance.
\end{abstract}

$\textbf{Keywords}$: Derivative-free optimization; Linear interpolation; Quadratic interpolation; Minimum Frobenius norm; Dynamic model selection

$\textbf{AMS Subject Classification}$: primary  90C56; secondary  49M25, 65K05.



\section{Introduction}\label{introduction}

Derivative-Free optimization (DFO) focuses on designing methods to solve optimization problems without evaluating derivatives or gradients of the objective function. This is particularly applicable when the objective function is a black-box function or an oracle function, so the only available information is the value of the objective function for an input point. The study of DFO methods has grown in recent years.  This is partly due to the flexibility of DFO methods across a variety of applied problems \cite{booker-etal-1998,marsden-feinstein-taylor-2008,marsden-wang-dennis-moin-2007,hare-2010,Bigdeli-Hare-Tesfamariam-2012-SPIE} (among many other examples) and partly due to the development of mathematics to ensure convergence \cite{Audet-Dennis-2002,custodio2008,C.Bogani2009,W.L.Hare2014} (among many other examples).  For a general overview of the DFO methods, along with a comprehensive study of many convergence results, see \cite{conn-scheinberg-vicente-2009} and the many references therein.

There are two main families of DFO methods: model-based methods and direct search methods. Direct search methods, at each iteration, sample the objective function at a finite number of points and act based on those function values without any derivative approximation.  A wide range variants exist, and many enjoy strong theoretical convergence analysis.  However, these are not the focus of this paper.

Model-based DFO methods use past information about the objective function $f$ to create a model function $\tilde{f}$ that approximates $f$.  Recent years have seen a significant amount of research focused on how to construct a `good' model function for use in DFO \cite{conn-scheinberg-vicente-2008,conn-scheinberg-vicente-2008b,ScheinbergToint2010,Gumma2014} (among others).

Among the most common used in practice are models formed from linear interpolation \cite{Powell1992}, quadratic interpolation \cite{ConnScheinbergToint1997,Powell2002}, or minimum Frobenius norms \cite{Powell2004,CustodioRochaVicente2010}.  (Although, it should be noted that other techniques exist, for example radial basis function models \cite{WildRegisShoemaker2008,YaoChenHuangvanTooren2014}, stochastic models \cite{BandeiraScheinbergVicente2014}, and models based on the Gaussian process \cite{ScottFrazierPowell2011}.) In $\R^n$, linear interpolation uses $n+1$ (well-poised) points to create an approximate gradient \cite[\S 2.3]{conn-scheinberg-vicente-2009}.  Alternatively, quadratic interpolation requires $(n+1)(n+2)/2$ (well-poised) points, but creates an approximate gradient and an approximate Hessian \cite[\S 3.4]{conn-scheinberg-vicente-2009}.   Details on linear interpolation and quadratic interpolation are given in Section \ref{models}.

As linear interpolation only provides approximate gradients, one would expect that, in terms of the number of iterations the resulting algorithm will probably converge similar to typical first order methods (i.e., linearly).  Conversely, as quadratic interpolation provides the benefit of an approximate Hessian, one might conjecture that the resulting algorithm converges similar to a second order method (i.e., quadratically).  Of course, in DFO, optimizers typically measure convergence in terms of number of function calls, not in terms of iterations.  (The assumption is that a function call takes significantly longer than any other portion of the algorithm.)  This means, one can take $(n+2)/2$ iterations of a method using linear interpolation for every single iteration of a method using quadratic interpolation.

This led to the development and use of minimum Frobenius norm models \cite{Powell2004,CustodioRochaVicente2010}.  Minimum Frobenius norm models provide some balance between the extremes of linear and quadratic interpolation.  Using between $n+1$ and $(n+1)(n+2)/2$ (well-poised) points, a minimum Frobenius norm model creates an approximate gradient that is more accurate than linear interpolation, along with an approximate Hessian that is less accurate than quadratic interpolation. (Details on minimum Frobenius norms are given in Section \ref{models}.) It is hoped that this balance can lead to DFO algorithms with strong convergence rates in terms of number of iterations, without the need for excessive function calls per iteration.

The availability of these three common models raises some natural questions in model-based DFO research.  First, how many points should be used to create the model function, i.e., how many points should be in the {\em sample set}?  Second, should the number of points in the sample set be static throughout the algorithm or can it be dynamically updated based on how the algorithm performed in the previous iteration?  In this paper, we present the results of a case study examining these questions.  Our case study focuses on the Derivative-Free Proximal Point (\DFPP) algorithm introduced by Hare and Lucet in 2014 \cite{W.L.Hare2014}.  We test 64 basic strategies for updating the size of the sample set, and compare the strategies across 60 test problems.  Numerical results suggest that the number of points in the sample set has a huge impact on algorithm performance, and adaptive strategies can both help and harm algorithm convergence.

The remainder of this paper is organized as follows.  In Section \ref{models}, we provide background details on how to construct linear interpolation, quadratic interpolation, and minimum Frobenius norm models.  In Section \ref{dfpp}, we outline the \DFPP~algorithm and present our adaptive strategies for determining the size of the interpolation set.  In Section \ref{numerics}, we discuss our numeric tests, present the results, and provide some qualitative remarks.  In Section \ref{conclusion}, we provide some conclusions.  Tables of results appear in the appendix.

Throughout, $B_{\Delta}(y^0)$ denotes the closed ball of radius $\Delta$ centred at $y^0$: $B_{\Delta}(y^0)=\{x:\|x-y^0\|\leq\Delta\}$,

\section{Model Construction Techniques}\label{models}
Let $\mathcal{P}_n^d$ be the space of polynomials of degree less than or equal to $d$ in $\mathbb{R}^n$.  For $d=1$ and $d=2$, the dimension of this space  is $\mathrm{dim}(\mathcal{P}_n^1)=n+1$ and $\mathrm{dim}(\mathcal{P}_n^2)=\frac{(n+1)(n+2)}{2}$, respectively. A basis $\Phi=\{\phi_1,...,\phi_q\}$ of $\mathcal{P}_n^d$ is a set of $q$ polynomials of degree less than or equal to $d$ such that $q=\mathrm{dim}(\mathcal{P}_n^d)$ and the polynomials span $\mathcal{P}_n^d$. If $\Phi$ is a basis in  $\mathcal{P}_n^d$, then any polynomial $m \in\mathcal{P}_n^d$ can be formulated as $m(x)=\sum_{j=1}^q \alpha_j\phi_j(x)$, where $\{\alpha_1, \alpha_2, ... \alpha_q\}$ is a uniquely determined set of real coefficients.

We say the polynomial $m$ interpolates the function $f$ at a given point $y$ if $m(y)=f(y)$. Suppose we are given a set $Y=\{y^0,y^1,...,y^p\}$ of interpolation points, and we seek a polynomial, $m$, with degree less than or equal to $d$ that interpolates a given function $f$ at the points in $Y$. Since it must be possible to write $m$ in the form of $\sum_{j=1}^q \alpha_j\phi_j$, we seek interpolation coefficients, $\alpha_j$, that satisfy the interpolation conditions
\begin{align}\label{eq:interpolate}
m(y^i)=\sum_{j=1}^q \alpha_j\phi_j(y^i)=f(y^i), i=0,...,p.
\end{align}

Conditions \eqref{eq:interpolate} form a linear system in terms of the interpolation coefficients, which we will write in matrix form as $M(\Phi,Y)\alpha_{\phi}=f(Y),$
where
    $$\begin{array}{c}M(\Phi,Y)=\left(\begin{array}{cccc}
                             \phi_0(y^0) & \phi_1(y^0) & \ldots & \phi_q(y^0) \\
                             \phi_0(y^1) & \phi_1(y^1) &\ldots  & \phi_q(y^1) \\
                             \vdots& \vdots & \vdots & \vdots \\
                             \phi_0(y^p) & \phi_1(y^p) & \ldots & \phi_q(y^p) \\
                \end{array} \right),
    \\
    \alpha_{\phi}=\left(\begin{array}{c}
                   \alpha_0 \\ \alpha_1 \\ \vdots \\ \alpha_q \\
                \end{array} \right),
\quad \mbox{and} \quad f(Y)=\left(\begin{array}{c}
                         f(y^0) \\ f(y^1) \\ \vdots \\ f(y^p) \\
                \end{array} \right).
    \end{array}$$
If conditions \eqref{eq:interpolate} have a unique solution, then their solution provides $m$.  If conditions \eqref{eq:interpolate} have multiple solutions, then $m$ is said to be under determined, and we must impose additional  conditions to select $m$.  If conditions \eqref{eq:interpolate} have no solution, then $m$ cannot exist, but a least-squares solution could be used to find an approximate solution and create an $m$ that approximates $f$ in the sense of statistical regression.

\subsection{Linear Interpolation}
Linear interpolation sets the maximum degree of the polynomial to $d=1$; i.e., linear interpolation applies in $\mathcal{P}_n^1$.  The natural basis for this space is $\Phi = \{1, x_1, x_2, ...x_n\}$.  Our interpolation conditions (equation \eqref{eq:interpolate}), can be simplified to
$$M(Y)\alpha = f(Y)$$
where
$$ M(Y) = \left(\begin{array}{cccc}
                             1 & y^0_1 & \ldots & y^0_n \\
                             1 & y^1_1 & \ldots & y^1_n \\
                             \vdots& \vdots & \vdots & \vdots \\
                             1 & y^n_1 & \ldots & y^n_n \\
                \end{array}\right),
\quad
    \alpha = \left( \begin{array}{c}
                   \alpha_0 \\\alpha_1 \\\vdots \\\alpha_n \\
             \end{array}\right),
\quad \mbox{and} \quad f(Y)=\left(
             \begin{array}{c}
                   f(y^0) \\f(y^1) \\\vdots \\f(y^n) \\
            \end{array}\right).$$
Clearly, $M(Y)$ is invertible if and only if conditions \eqref{eq:interpolate} have a unique solution.  More importantly, the error in the gradient approximation from linear interpolation can be quantified in terms of several constants and the approximate radius of the sample set.

\begin{Theorem}\textbf{(Error bound for Linear interpolation)}\cite[Thm 2.11]{conn-scheinberg-vicente-2009}\label{thm:linint}

Let $Y =\{y^0, y^1,\ldots, y^n\}\subseteq \mathcal{R}^n$ be poised for linear interpolation.  Define $\Delta=\Delta (Y)=\max_{1\leq i\leq n} \| y^i-y^0\|$. Suppose the function $f$ is continuously differentiable in an open domain $\Omega$ containing $B_{\Delta}(y^0)$, and $\nabla f$ is Lipschitz continuous in $\Omega$ with constant $\nu>0$. Let $m$ be the linear function that interpolates $f$ over all points in $Y$.  Then for all points $y \in B_{\Delta}(y^0)$, we have
    $$\|\nabla f(y)-\nabla m(y)\|\leq \kappa_{eg}\Delta,$$
where $\kappa_{eg}$ is a constant based on $\nu$, $n$, and the geometry of the interpolation set.
\end{Theorem}

\subsection{Quadratic Interpolation}

Quadratic interpolation sets the maximum degree of the polynomial to $d=2$; i.e., quadratic interpolation applies in $\mathcal{P}_n^2$.  One natural basis for quadratic interpolation is
    \[\Phi = \{1, x_1, x_2, ...x_n, \frac{1}{2}(x_1)^2, x_1 x_2, ... x_1 x_n, \frac{1}{2}(x_2)^2, x_2 x_3, ... \frac{1}{2}(x_n)^2\}.\]
Using this basis, one can again write conditions \eqref{eq:interpolate} as a linear system. (For the sake of space, we do not rewrite the system here.)  The system results in the matrix $M(\Phi, Y)$ being a $(n+1)(n+2)/2 \times (n+1)(n+2)/2$ square matrix.  Like linear interpolation, the error in the gradient approximation from quadratic interpolation can be quantified using several constants and the approximate radius of the sample set.

\begin{Theorem}\textbf{(Error bounds for quadratic interpolation)}\cite[Thm 3.16]{conn-scheinberg-vicente-2009}\label{thm:quadint}

Let $Y =\{y^0, y^1,\ldots, y^p\}\subseteq \mathcal{R}^n$ be poised for quadratic interpolation.  Define $\Delta=\Delta (Y)=\max_{1\leq i\leq p} \| y^i-y^0\|$. Suppose the function $f$ is twice continuously differentiable in an open domain $\Omega$ containing $B_{\Delta}(y^0)$, and $\nabla^2 f$ is Lipschitz continuous in $\Omega$ with constant $\nu_2>0$.  Let $m$ be the quadratic function that interpolates $f$ over all points in $Y$.  Then, for all points $y \in B_{\Delta}(y^0)$, we have
    \begin{align*}
    \|\nabla^2 f(y)-\nabla^2 m(y)\|&\leq \kappa_{eh}\Delta, ~\mbox{and}\\
    \|\nabla f(y)-\nabla m(y)\|&\leq \kappa_{eg}\Delta^2,\\
    \end{align*}
where $\kappa_{eh}$ and $\kappa_{eg}$ are constants based on $\nu_2$, $p$, and the geometry of the interpolation set.\end{Theorem}

\subsection{Minimum Frobenius Norm Models}

Minimum Frobenius norm models set the maximum degree of the polynomial to $d=2$, but work in the case when the number of interpolation points is less than the $(n+1)(n+2)/2$ required for quadratic interpolation.

In this case, the $M(\Phi,Y)$ defining the interpolating conditions has more columns than rows and the interpolation polynomials are no longer unique.

Let us split the natural basis $\Phi$ into linear and quadratic parts: ${\Phi}_L=\{1,x_1,...,x_n\}$ and ${\Phi}_Q=\{\frac{1}{2}x_1^2, x_1 x_2,..., \frac{1}{2}x_n^2\}.$  The interpolation model can now be written as
    $$m(x)=\alpha_L^T{\Phi}_L(x)+\alpha_Q^T{\Phi}_Q(x),$$
where $\alpha_L$ and $\alpha_Q$ are the appropriate parts of the coefficient vector $\alpha$.

In a DFO framework with under determined interpolation, it is desirable to construct accurate linear models and then enhance them with curvature information, hoping that the actual accuracy of the model is better than that of a purely linear model. (Hence, it is important to construct sample sets that are poised for linear interpolation.)

Since the interpolation set is too small to create a unique quadratic interpolation, we must impose some additional requirements to determine the final model.  As our Hessian approximation will be of a lower accuracy than our gradient approximation, in derivative-free optimization it makes sense to seek a model for which the norm of the Hessian is small or moderate. Therefore, we define the minimum Frobenius norm solution as a solution to the following optimization problem in $\alpha_L$ and $\alpha_Q$:
    $$\min \frac{1}{2}\parallel\alpha_Q\parallel_2^2$$
    \begin{equation}
    M({\Phi}_L,Y)\alpha_L+M({\Phi}_Q,Y)\alpha_Q=f(Y).
    \end{equation}
The name minimum Frobenius norm solution comes from the equivalence of minimizing the norm of $\alpha_Q$ and minimizing the Frobenius norm of the Hessian of $m$.

The condition for the existence and uniqueness of the minimum Frobenius norm model is that the following matrix is nonsingular
\begin{align}
F({\Phi},Y)=\left(
  \begin{array}{cc}
    M({\Phi}_Q,Y)M({\Phi}_Q,Y)^T & M({\Phi}_L,Y) \\
    M({\Phi}_L,Y)^T & 0\\
  \end{array}
\right).
\end{align}
We say that a set $Y$ is poised for minimum Frobenius norm interpolation if problem (2) has a unique solution or, equivalently, if the matrix $F({\Phi},Y)$ is non-singular. Like linear and quadratic interpolation, the error bounds of the resulting model can be bounded using constants and the approximate radius of the sample set.

\begin{Theorem}\textbf{(Error bounds for minimum Frobenius norm models)}\cite[Thm 5.4]{conn-scheinberg-vicente-2009}\label{thm:Frobint}

Let $Y =\{y^0, y^1,\ldots, y^p\}\subseteq \mathcal{R}^n$ be poised for minimum Frobenius norm interpolation.  Define $\Delta=\Delta (Y)=\max_{1\leq i\leq p} \parallel y^i-y^0\parallel$. Suppose the function $f$ is continuously differentiable in an open domain $\Omega$ containing $B_{\Delta}(y^0)$, and $\nabla f$ is Lipschitz continuous in $\Omega$ with constant $\nu>0$.  Let $m$ be the quadratic function that results from the minimum Frobenius norm interpolation of $f$ over all points in $Y$. Then, for all points $y \in B_{\Delta}(y^0)$, we have
    $$
    \|\nabla f(y)-\nabla m(y)\| \leq \kappa_{eg}\Delta,
    $$
and $\kappa_{eg}$ is a constants based on $\nu$, $p$, the geometry of the interpolation set, and the norm of the model Hessian.
\end{Theorem}

\section{Derivative-Free Proximal Point Method}\label{dfpp}

Much like Newton's method is a standard tool for solving smooth optimization problems, proximal point algorithms can be viewed as an analogous tool for nonsmooth optimization. The basic (theoretical) method solves the minimization of $f$ through iterative solutions to the proximal point problem
    \begin{equation}\label{eq:PPdef}
    x^{k+1} = \mathrm{prox}_r f(x^k) ~~\mbox{where}~ \mathrm{prox}_r f(x^k) = \arg\min\{f(y) + \frac{r}{2}\|y - x^k\|^2\}.
    \end{equation}
In practice, it is unnecessary to solve $\mathrm{prox}_r f(x^k)$ exactly, which has lead to a variety of practical implementations based on the proximal point framework  \cite{Hiriart-Urruty-Lemarechal-1993a,Makela-2002,OliveiraSagastizabal2014} (and references therein).

Most common are the proximal-bundle methods, where the objective function $f$ is replaced by a sequence of piecewise linear model functions $f_k$, see \cite{Hiriart-Urruty-Lemarechal-1993a,Makela-2002,OliveiraSagastizabal2014}. Such methods essentially replace the minimization of $f$ with a sequence of quadratic programming problems.

Recently, Hare and Lucet introduced a Derivative-Free Proximal Point (\DFPP) method \cite{W.L.Hare2014}. Within the framework, $x^k$ denotes the prox-centre of the algorithm during iteration $k$. At each iteration, the algorithm shall make use of a sample set $Y=\{y^0, y^1,...,y^p\}\subseteq \R^n$ with $y^0 = x^k$ to construct a model of the objective function.  In the algorithm, the approximate sampling radius of $Y$ is defined $\Delta(Y)=\max_{y^i\in Y}\| y^i-y^0\|$, and $\lambda_n(H)$ is used to denote the minimum eigenvalue of $H$.  Pseudo-code of the \DFPP~algorithm follows.

\textbf{Derivative-Free Proximal Point Method (DFPP)}
\begin{enumerate}
	\item[0.] {\sc Initialize}: Set $k=0$ and input
		\begin{tightenumerate}
			\item[] $x^0$ - an initial prox-centre,
			\item[] $Y^0$ - an initial poised interpolation set with $x^0 \in Y^0$,
			\item[] $r^0$ - an initial prox-parameter, $r^0 > 0$
			\item[] $m$ - an Armijo-like parameter, $0 < m < 1$,
			\item[] $\Gamma$ - a minimal radius decrease parameter, $0 < \Gamma < 1$,
			\item[] $\Rtol$ - stopping tolerance for prox-parameter, $\Rtol > 0$,
			\item[] $\Dtol$ - stopping tolerance for search radius, $\Dtol \geq 0$, and
			\item[] $\Gtol^\nabla, \, \Gtol^\Delta$ - stopping tolerances for approximate gradient, $\Gtol^\nabla, \, \Gtol^\Delta > 0$.
		\end{tightenumerate}
	\item[1.] {\sc Model and Stopping Conditions}: \\
		Create $q^k$, a model of $f$ over $Y^k$:
			\[q^k(x) := a^k + \langle g^k, x \rangle + \frac{1}{2}\langle x, H^k x\rangle .\]
		If $\|\nabla q^k(x^k)\| < \Gtol^\nabla$ and $\Delta(Y^k) < \Gtol^\Delta$, then STOP (`success').\newline
		If $\Delta(Y^k) < \Dtol$, then STOP.
	\item[2.] {\sc Prox-feasiblity Check}: \\
		If  $r^k \leq -\lambda_n(H^k)$, then ($q^k + r^k \frac{1}{2}\|\cdot\|^2$ is not strictly convex):
			\begin{tightenumerate}
			\item reset $r^k = -\lambda_n(H^k)+1$,
			\end{tightenumerate}
	\item[3.] {\sc Prox Trial Point}: \\
		Compute the trial point
			\[\{\tilde{x}^k\} = \mathrm{prox}_{r^k} q^k (x^k) = \{(H^k + r^k \mbox{Id})^{-1}(r^k x^k - g^k)\}.\]
		Compute the predicted decrease
			\[\delta^k = q^k(x^k) - q^k(\tilde{x}^k).\]
	\item[4.] {\sc Serious/Null Check}:\\
		If $f(\tilde{x}^k) \leq f(x^k) - m \delta^k$, then declare a {\em serious step}:
			\begin{tightenumerate}
            \item select $x^{k+1}$ such that $f(x^{k+1}) \leq f(x^k) - m \delta^k$,
			\item generate an interpolation set $Y^{k+1}$ such that $x^{k+1}\in Y^{k+1}$ and $\Delta(Y^{k+1}) \leq \Delta(Y^k)$.
			\end{tightenumerate}
		Else (if $f(\tilde{x}^k) > f(x^k) - m \delta^k$), then declare a {\em null step}:
			\begin{tightenumerate}
			\item if $\tilde{x}^k \notin B_{\Delta(Y^k)} (x^k)$, then declare the null step to be {\em type 1},  increase $r^{k+1} \rightarrow 2 r^k$, set $x^{k+1} = x^k$ and $Y^{k+1}$ with $Y^k \subseteq Y^{k+1}$ and $\Delta(Y^{k+1})=\Delta(Y^k)$, \\
			\item if $\tilde{x}^k \in B_{\Delta(Y^k)} (x^k)$,  then declare the null step to be {\em type 2}, set $x^{k+1} = x^k$ and generate an interpolation set $Y^{k+1}$ such that $x^{k+1}\in Y^{k+1}$ and $\Delta(Y^{k+1}) \leq \Gamma \Delta(Y^k)$.
			\end{tightenumerate}
	\item[5.] {\sc Loop}:

	Increment $k \rightarrow k+1$ and return to Step 1.
	\end{enumerate}

\subsection{Adaptive Strategies in DFPP}\label{strategies}

The \DFPP~framework has two interesting features that make it well-suited to exploring adaptive updating of the number of points in the sample set at each iteration.

First, in order for the algorithm to converge the model functions $q^k$ must satisfy the following assumption (see \cite{W.L.Hare2014}).

\begin{assumption}\label{assume:linear}
Assume $f  \in \mathcal{C}^1$.  Furthermore, assume that there exists constants $C$ and $M$ such that, for any point $y^0$ and any sampling radius $\Delta > 0$, we are able to generate a sampling set $Y = \{ y^0, y^1, \ldots y^p\} \subseteq \R^n$ and a corresponding quadratic model function $q$ such that $\Delta(Y) = \Delta$ and
	\[ \begin{array}{rcll}
	\| f(y) - q(y) \| &\leq& C \Delta^2 & \mbox{for all}~ y \in B_\Delta(y^0), \\
	\| \nabla f(y) - \nabla q(y) \| &\leq& C  \Delta & \mbox{for all}~ y \in B_\Delta(y^0), ~\mbox{and} \\
	\|\nabla^2 q(y) \| &\leq& M.
	\end{array}\]
\end{assumption}

\noindent Note that the Assumption \ref{assume:linear} can be satisfied through
    \begin{itemize}
    \item[a)] a linear interpolation model, by noting $\nabla^2 q = 0$ in this case,
    \item[b)] a quadratic interpolation model, see \cite[Lem 3.1]{W.L.Hare2014}, or
    \item[c)] a minimum Frobenius norm model, see \cite[p. 209]{W.L.Hare2014}.
    \end{itemize}
Moreover, the algorithm can use a different model construction technique at each iteration, without compromising the convergence analysis.

The second interesting feature of the \DFPP~framework is that the algorithm ends step 4 with one of three possible declarations: {\em serious step, null step type 1,} or {\em null step type 2}.  In a {\em serious step}, the algorithm was successful in finding a new proximal centre, which shows notable decrease over the previous proximal centre.  In a {\em null step type 1}, the algorithm was unable to find a new proximal centre and the predicted new centre was outside of the radius of accuracy for the model.  In this case, the prox-parameter is increased, but the old sample set can be reused (possibly with additional points added).  Finally, in a {\em null step type 2}, the algorithm was unable to find a new proximal centre, despite the fact that the predicted new centre was inside of the radius of accuracy for the model.  In this case, the desired accuracy of the model is increased (i.e., $\Delta(Y)$ is decreased), and a new model must be constructed.

As each outcome suggests a different situation, and each iteration can use a different model construction technique, the \DFPP~algorithm naturally lends itself to the idea of using an adaptive strategy for selecting sample set size at each iteration.  The adaptive strategies explored in this paper can be viewed as selecting the number of points in the sample set in the next iteration based on the declaration in step 4:
   \begin{enumerate}
   \item[i.] if step 4 of the \DFPP~method declares a serious step, then use a sample set of size $N_{\tt s}$ in the next iteration,
   \item[ii.] if step 4 of the \DFPP~method declares a null step type 1, then use a sample set of size $N_{\tt n1}$ in the next iteration, and
   \item[iii.] if step 4 of the \DFPP~method declares a null step type 2, then use a sample set of size $N_{\tt n2}$ in the next iteration.
   \end{enumerate}
For our testing, $N_{\tt s}, N_{\tt n1},$ and $N_{\tt n2}$ are each taken from
    \[\left\{n+1, 2n+1, \left\lfloor\frac{1}{2} \left((n+1)+\frac{(n+1)(n+2)}{2}\right) \right\rfloor, \frac{(n+1)(n+2)}{2}\right\},\]
where $n$ is the problem dimension.  In order to simplify presentation, we use the following notation for these four strategies
    \[\begin{array}{rcl}
    \lin &=& n+1, \\
    \bilin &=& 2n + 1,\\
    \mean &=& \left\lfloor\frac{1}{2} \left((n+1)+\frac{(n+1)(n+2)}{2}\right) \right\rfloor,\\
    \qua &=& \frac{(n+1)(n+2)}{2}.\\
    \end{array}\]
As there are 3 possible ways to conclude step 4, and we examine 4 different strategies for each conclusion, we explore a total of $4^3 = 64$ different adaptive strategies.

\subsection{Sample Set Construction Techniques}\label{sampleset}

The error bounds provided in Theorems \ref{thm:linint}, \ref{thm:quadint}, and \ref{thm:Frobint}, all rely on the ``geometry of the interpolation set''.  This phrase actually hides a deep literature on the topic.  While, in order for the sample set to be poised, we require an appropriate matrix, $M(\Phi, Y)$ or $F(\Phi, Y)$, to be invertible, in practice it is important that this matrix is `stable'. This stability is dependent on the ``geometry of the interpolation set''.  Details on quantifying the geometry of the interpolation set, and how to control the quality of this geometry, are outside of the scope of this work (we refer interested readers to  \cite[Chpt 2--6]{conn-scheinberg-vicente-2009}).  Nonetheless, some comments are in order.

For the numerical testing in this paper, we use Algorithms 6.2 and 6.3 of \cite{conn-scheinberg-vicente-2009} to construct our interpolation sets and improve their geometry.  The interpolation set $Y^k$ is built based on three possible declarations in step 4.
\begin{enumerate}
    \item[i.] If step 4 of the \DFPP~method declares a serious step, then
        \begin{tightenumerate}
            \item place $x^{k+1}$ into an interpolation set $Y$ in position of $x = y_0$,
            \item place any previously sampled points in $B_{\Delta(Y^k)}(x^{k+1})$ into $Y$,
            \item use \cite[Alg 6.2 \& 6.3]{conn-scheinberg-vicente-2009} to expand $Y$ into a well-poised set of $(n+1)(n+2)/2$ points, and
            \item select $N_{\tt s}$ points from $Y$ to define $Y^{k+1}$.
        \end{tightenumerate}

    \item[ii.] If step 4 of the \DFPP~method declares a null step type 1, then
        \begin{tightenumerate}
        \item if the number of points in $Y^k$ is greater or equal to $N_{\tt n1}$, then select $N_{\tt n1}$ points from $Y^k$ to define $Y^{k+1}$,
        \item if the number of points in $Y^k$ is less than $N_{\tt n1}$, then use \cite[Alg 6.2 \& 6.3]{conn-scheinberg-vicente-2009} to create a well-poised set of $(n+1)(n+2)/2$ points and select $N_{\tt n1}$ points from it to define $Y^{k+1}$.
        \end{tightenumerate}

    \item[iii.] If step 4 of the \DFPP~method declares a null step type 2, then
        \begin{tightenumerate}
            \item place $x^{k+1}$ into an interpolation set $Y$ in position of $x = y_0$,
            \item place any previously sampled points in $B_{\Gamma\Delta(Y^k)}(x^{k+1})$ into $Y$,
            \item use \cite[Alg 6.2 \& 6.3]{conn-scheinberg-vicente-2009} to expand $Y$ into a well-poised set of $(n+1)(n+2)/2$ points, and
            \item select $N_{\tt n2}$ points from $Y$ to define $Y^{k+1}$.
        \end{tightenumerate}
\end{enumerate}

In all of the cases above, the selection of the final subset of $N_{\tt x}$ points from $Y$ is done by taking the first $N_{\tt x}$ points, and then a safety check is used to ensure the final interpolation set is well-poised.  If it is not, then we use \cite[Alg 6.2 \& 6.3]{conn-scheinberg-vicente-2009} to create a well-poised set of $(n+1)(n+2)/2$ points and select the first $N_{\tt x}$ points from it to define $Y^{k+1}$ (again with a safety check to ensure well-poised).

\section{Numerical Results}\label{numerics}
The \DFPP~method is implemented in MATLAB \cite{W.L.Hare2014}.  Minor adaptations to the original code allowed for the adaptive strategies in Subsection \ref{strategies} to be incorporated. Minor tuning to select algorithmic parameters was performed.  Specifically, Armijo-like parameters $m\in \{0.1, 0.5, 0.9\}$ were tested and radius decrease parameters $\Gamma \in \{0.25, 0.5\}$ were tested.  The values $m=0.1$ and $\Gamma = 0.5$ provided the best overall performance across all strategies. The initial prox-parameter was set to $r^0=1$. (As improvement based on these parameters was extremely minor, we do not present results from other parameter combinations; however, these results are available through contacting the corresponding author.) Finally, in Step 4, if a serious step is declared, the user has the option of performing a line search (or other search method) to seek $x^{k+1}$ that provides some further improvement over $\tilde{x}^k$.  We tested using no additional searching and using a backtracking line search.

The strategies were tested on the 60 problems from \cite{AliKhomZab2005,CODE,More-Garbow-Hillstrom-1981}.  Test problems were separated into two groups: low dimension and high dimension. Table \ref{tab:testproblems} lists the name and the dimension of each test problem.

\begin{table}[ht]\caption{Test problems.}\label{tab:testproblems}
\begin{center}
\begin{tabular}{|l|c|p{1.25cm}|l|c|}\cline{1-2}\cline{4-5}
\multicolumn{2}{|c|}{Low Dim. Problems} && \multicolumn{2}{|c|}{High Dim. Problems} \\
\cline{1-2}\cline{4-5}
Function Name	&	Dim. &&	Function Name & Dim. \\\cline{1-2}\cline{4-5}
Bard & 3 && Ackley & 10 \\\cline{1-2}\cline{4-5}
Beale & 2 && Ackley & 20 \\\cline{1-2}\cline{4-5}
Biggs EXP6 & 6 && Arrowhead & 10 \\\cline{1-2}\cline{4-5}
Box 3D & 3 && Arrowhead & 20 \\\cline{1-2}\cline{4-5}
Brown almost-linear & 3 &&  Epistatic Michalewicz & 10 \\\cline{1-2}\cline{4-5}
Brown \& Dennis & 4 && Exponential & 10 \\\cline{1-2}\cline{4-5}
Brown badly scaled & 2 && Exponential & 20 \\\cline{1-2}\cline{4-5}
Broyden banded & 3 && Griewank & 10 \\\cline{1-2}\cline{4-5}
Broyden tridiagonal & 3 && Griewank & 20 \\\cline{1-2}\cline{4-5}
Discrete boundary value & 3 && Levy \& Montalvo I & 10 \\\cline{1-2}\cline{4-5}
Discrete integral eq. & 3 && Levy \& Montalvo I & 20 \\\cline{1-2}\cline{4-5}
Freudenstein \& Roth & 2 && Levy \& Montalvo II & 10 \\\cline{1-2}\cline{4-5}
Gaussian & 3 && Levy \& Montalvo II & 20 \\\cline{1-2}\cline{4-5}
Gulf & 3 && Modified Langerman & 10  \\\cline{1-2}\cline{4-5}
Helical valley & 3 && Neumaier 3 & 10  \\\cline{1-2}\cline{4-5}
Jennrich \& Sampson & 2 && Neumaier 3& 20 \\\cline{1-2}\cline{4-5}
Kowalik \& Osborne & 4 && Powell singular & 12 \\\cline{1-2}\cline{4-5}
Linear rank-1 & 3 && Powell singular & 20 \\\cline{1-2}\cline{4-5}
Linear rank full & 4 && Rastrigin & 10 \\\cline{1-2}\cline{4-5}
Linear rank-1 with zeros & 3 && Rastrigin & 20 \\\cline{1-2}\cline{4-5}
Meyer & 3 && Rosenbrock & 10 \\\cline{1-2}\cline{4-5}
Osborne I & 5 && Rosenbrock & 20 \\\cline{1-2}\cline{4-5}
Penalty I & 4 && Sinusoidal & 10 \\\cline{1-2}\cline{4-5}
Penalty II & 4 && Sinusoidal & 20 \\\cline{1-2}\cline{4-5}
Powell badly scaled & 2 &&  Variably dimensional & 10 \\\cline{1-2}\cline{4-5}
Rosenbrock & 2 &&  Variably dimensional & 20 \\\cline{1-2}\cline{4-5}
Trigonometric & 3 && Wood & 10 \\\cline{1-2}\cline{4-5}
Variably dimensional & 3 && Wood & 20 \\\cline{1-2}\cline{4-5}
Watson & 6 && Zakharov & 10 \\\cline{1-2}\cline{4-5}
Wood & 4 && Zakharov & 20 \\\cline{1-2}\cline{4-5}
\end{tabular}
\end{center}\end{table}

For each test problem, each strategy was run until a total of $100 n$ function calls was exceeded, where $n$ is the dimension of the test problem.

In order to rank the strategies, we consider the following {\em improvement} metric
    \[\imp(N_{\tt s}, N_{\tt n1}, N_{\tt n2}) = \sum_{\mathcal{P}} \min\left\{ -\log_{10}\frac{|(f-fbest)|}{|(f_0-fbest)|}, 16 \right\}\]
where $\mathcal{P}$ is the set of all test problems, $f$ is the objective function value obtained by \DFPP, $fbest$ is the best known objective function value, and $f_0$ is the initial objective function value.  The value $ -\log_{10}\frac{|(f-fbest)|}{|(f_0-fbest)|}$ can loosely be interpreted as the number of new digits of accuracy (in function value) obtained on a given test problem.  The minimization with $16$ deals with the (few) problems that end up being solved exactly and return unrealistic values like $-\log_{10}\frac{|(f-fbest)|}{|(f_0-fbest)|} > 1000$.  Without capping, these problems can massively skew the data analysis.  Finally, these values are summed over all test problems, to give each strategy a total improvement.

The aggregate results when no line search was used appear in Tables \ref{tab:nolineall} to \ref{tab:nolinehigh} (in the Appendix) and the aggregate results when a backtracking line search was used appear in Tables \ref{tab:withlineall} to \ref{tab:withlinehigh} (in the Appendix)\footnote{For more detailed results please contact the corresponding author.}.  Tables \ref{tab:nolineall} and \ref{tab:withlineall} provide the results when all test problems are considered.  Tables \ref{tab:nolinelow} and \ref{tab:withlinelow} provide the results when only low dimension test problems are considered. Tables \ref{tab:nolinehigh} and \ref{tab:withlinehigh} provide the results when only high dimension test problems are considered.

\subsection{Interpretation of the Results}

To ease interpretation, Tables \ref{tab:nolineall} to \ref{tab:withlinehigh} are sorted from the highest $\imp$ value to the lowest $\imp$ value.

Examining Tables \ref{tab:nolineall} and \ref{tab:withlineall} (which contain all test problems grouped), we note that the line search has a strong positive impact on the performance of the algorithm.  This was also noted in \cite{W.L.Hare2014}.  In fact, the line search is so effective that the worst result in  Table \ref{tab:withlineall}, would rank $7^{\mbox{th}}$ if it were placed in Table \ref{tab:nolineall}.

Across all tables we see a common trend of $N_{\tt s} = \lin$.  That is, if a serious step occurred, then the next model should be as simple to create as possible.  This makes sense, as serious steps correspond with success and movement of the prox-centre.  If a serious step occurred, then the next model essentially starts from scratch, so it makes sense to build a simple model and only increase complexity if the next iteration induces a null step.

Comparing low dimension to high dimension problems presents some enlightening results.  In both Tables \ref{tab:nolinelow} and  \ref{tab:withlinelow} (low dimensions), we see 4 of the top 5 strategies involve building complex models when a null step occurs (i.e., $N_{\tt n1} \in \{\mean, \qua\}$ or $N_{\tt n2} \in \{\mean, \qua\}$).  Conversely, in Tables \ref{tab:nolinehigh} and \ref{tab:withlinehigh} (high dimensions), we see complex models are generally avoided: in Table \ref{tab:nolinehigh}, $N_{\tt n1} \in \{\lin, \bilin\}$ and $N_{\tt n2} \in \{\lin, \bilin\}$ for all of the top 5, while in Table \ref{tab:withlinehigh}, $N_{\tt n1} \in \{\lin, \bilin\}$ and $N_{\tt n2} \in \{\lin, \bilin\}$ for 3 of the top 5.

\subsection{Data Profiles}

While Tables \ref{tab:nolineall} and \ref{tab:withlineall} provide some insight into the performance, they rely strongly on data aggregation.  As such, it is possible that certain problems are being solved to very high precision and skewing the results.  In this subsection we present data profiles \cite{MoreWild09} of select strategies.

Data profiles are designed to capture both speed and robustness of a solver, by plotting the portion of problems solved using less than or equal to $\alpha \times n+1$ function calls, where $\alpha$ is the number of `simplex gradients equivalents' used and $n+1$ represent the number of function calls required to create a simplex gradient in $\R^n$.  For further details on data profiles, we refer the reader to \cite{MoreWild09}.

\begin{figure}[ht]\begin{center}
\includegraphics[width=3in]{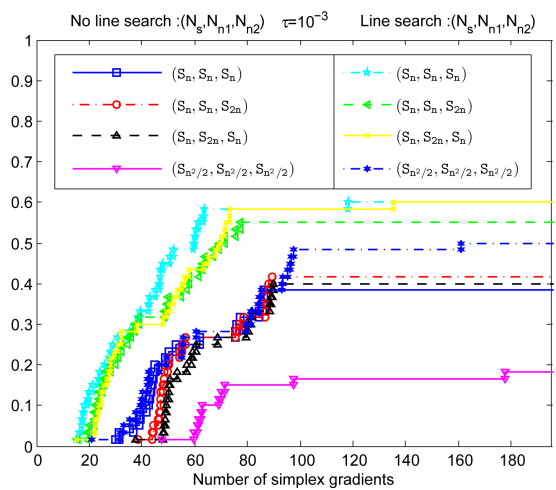}
\end{center}
\caption{Data profile for 8 strategies and a solving tolerance of $10^{-3}$.}\label{fig:data3}\end{figure}

\begin{figure}[ht]\begin{center}
\includegraphics[width=3in]{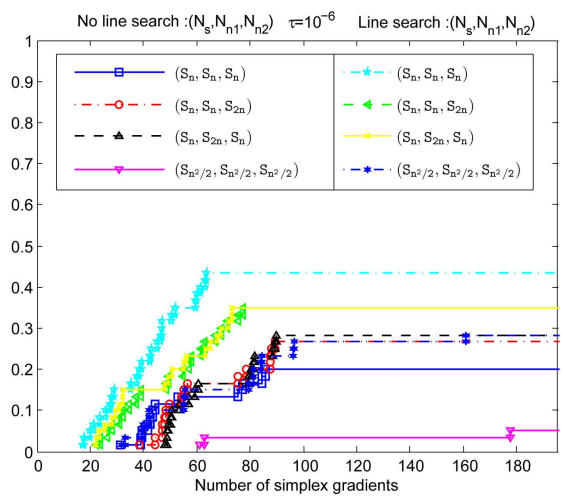}
\end{center}
\caption{Data profile for 8 strategies and a solving tolerance of $10^{-6}$.}\label{fig:data6}\end{figure}

With 128 strategies tested (64 adaptive approaches times 2 for the use/disuse of a line search), presenting all data profiles would result in an unreadable figure.  Instead, we present data profiles containing:
\begin{tightenumerate}
    \item the `most basic' non-adaptive strategy -- $N_{\tt s} = N_{\tt n1} = N_{\tt n2}=\lin$,
    \item the `most complex' non-adaptive strategy -- $N_{\tt s} = N_{\tt n1} = N_{\tt n2}=\qua$, and
    \item the top two adaptive strategies -- $N_{\tt s} = \lin$, $N_{\tt n1} = \lin$, $N_{\tt n2}= \bilin$, and $N_{\tt s} = \lin$, $N_{\tt n1} = \bilin$, $N_{\tt n2}= \lin$.
\end{tightenumerate}
Data profiles are created including the no line search and with a line search option.  The data profile for a solving tolerance of $10^{-3}$ appears in Figure \ref{fig:data3} and data profile for a solving tolerance of $10^{-6}$ appears in Figure \ref{fig:data6}.

Figure \ref{fig:data3} and \ref{fig:data6} show the expected results and a few surprising results.  First, it is again clear that the line search provides an excellent performance boost to the algorithm.  Examining just methods that use a line search, we note that strategy $N_{\tt s} = N_{\tt n1} = N_{\tt n2}=\lin$ outperforms all other methods.  However, examining the no line search methods, we see that for solving tolerance of $10^{-3}$ strategy $N_{\tt s} = N_{\tt n1} = \lin$, $N_{\tt n2}=\bilin$ outperforms the other no line search methods.  Meanwhile, for solving tolerance of $10^{-6}$ strategy $N_{\tt s} = N_{\tt n2} = \lin$, $N_{\tt n1}=\bilin$ outperforms the other no line search methods.  Neither of these victories are resounding, but it nonetheless suggests that adaptive strategies may have some place in future DFO algorithms.

\section{Conclusions}\label{conclusion}

Model-based DFO methods work by constructing local models of the objective function using a set of function evaluations.  In this paper, we explore the questions of how many function evaluations should be used to construct the model, and should this number be fixed, or adaptively selected by the algorithm?  We approach the question numerically, by making use of the flexibility and iteration decision structure within the \DFPP~algorithm of \cite{W.L.Hare2014}.  The results suggest that, for this algorithm, and this implementation, adaptive strategies can provide some improvement, particularly in higher dimensions.  However, the results also show that a poorly selected adaptive strategy can greatly hinder performance, both in low and high dimensions.  Finally, the results generally suggest that, for this algorithm and implementation, basic models using fewer function evaluations outperform complex models that require many function evaluations to build.

It should be noted that there are many model-building methods that were not considered in this paper: e.g., linear regression models, centered simplex gradients, radial basis functions, models based on the Gaussian process, etc.  Also, while past points within the sampling radius were used when building new models, advance techniques on minimum Frobenius norm based model updating was not applied within this paper.   This leaves significant opportunity for further research in this area.

\small

\subsection*{Acknowledgements}
WH research supported in part by the Natural Sciences and Engineering Research Council of Canada, Discovery Grant \#355571-2013.
\\
MJ expresses deep gratitude to Mathematics and Computer Science, Amirkabir University of Technology, for providing a scholarship and to the Centre for Optimization, Convex Analysis and Nonsmooth Analysis at the University of British Columbia for providing supportive and stimulating working environment during his stay as a Ph.D.\ Visiting Scholar.

\bibliography{DFPPM} 

\begin{thebibliography}{10}

\bibitem{AliKhomZab2005}
M.M. Ali, C.~Khompatraporn, and Z.B. Zabinsky.
\newblock A numerical evaluation of several stochastic algorithms on selected
  continuous global optimization test problems.
\newblock {\em Journal of Global Optimization}, 31(4):635--672, 2005.

\bibitem{Audet-Dennis-2002}
C.~Audet and J.E. Dennis, Jr.
\newblock Analysis of generalized pattern searches.
\newblock {\em SIAM J. Optim.}, 13(3):889--903 (electronic) (2003), 2002.

\bibitem{BandeiraScheinbergVicente2014}
A.S. Bandeira, K.~Scheinberg, and L.N. Vicente.
\newblock Convergence of trust-region methods based on probabilistic models.
\newblock {\em SIAM J. Optim.}, 24(3):1238--1264, 2014.

\bibitem{Bigdeli-Hare-Tesfamariam-2012-SPIE}
K.~Bigdeli, W.~Hare, and S.~Tesfamariam.
\newblock Optimal design of viscous damper connectors for adjacent structures
  using genetic algorithm and {N}elder-{M}ead algorithm.
\newblock In {\em Proceeding of SPIE conference on Smart Structures and
  Materials}, 2012.
\newblock Article 83410M (4 pages).

\bibitem{C.Bogani2009}
C.~Bogani, M.G. Gasparo, and A.~Papini.
\newblock Generating set search methods for piecewise smooth problems.
\newblock {\em SIAM J. Optim.}, 20(1):321--335, 2009.

\bibitem{booker-etal-1998}
A.J. Booker, J.E. Dennis, Jr., P.D. Frank, D.B. Serafini, and V.~Torczon.
\newblock Optimization using surrogate objectives on a helicopter test example.
\newblock In {\em Computational methods for optimal design and control
  ({A}rlington, {VA}, 1997)}, volume~24 of {\em Progr. Systems Control Theory},
  pages 49--58. Birkh\"auser Boston, Boston, MA, 1998.

\bibitem{ConnScheinbergToint1997}
A.~Conn, K.~Scheinberg, and P.~Toint.
\newblock On the convergence of derivative-free methods for unconstrained
  optimization.
\newblock In {\em Approximation theory and optimization ({C}ambridge, 1996)},
  pages 83--108. Cambridge Univ. Press, Cambridge, 1997.

\bibitem{conn-scheinberg-vicente-2008}
A.R. Conn, K.~Scheinberg, and L.N. Vicente.
\newblock Geometry of interpolation sets in derivative free optimization.
\newblock {\em Math. Program.}, 111(1-2, Ser. B):141--172, 2008.

\bibitem{conn-scheinberg-vicente-2008b}
A.R. Conn, K.~Scheinberg, and L.N. Vicente.
\newblock Geometry of sample sets in derivative-free optimization: polynomial
  regression and underdetermined interpolation.
\newblock {\em IMA J. Numer. Anal.}, 28(4):721--748, 2008.

\bibitem{conn-scheinberg-vicente-2009}
A.R. Conn, K.~Scheinberg, and L.N. Vicente.
\newblock {\em Introduction to derivative-free optimization}, volume~8 of {\em
  MPS/SIAM Series on Optimization}.
\newblock Society for Industrial and Applied Mathematics (SIAM), Philadelphia,
  PA, 2009.

\bibitem{custodio2008}
A.L. Cust{\'o}dio, J.E. Dennis, Jr., and L.N. Vicente.
\newblock Using simplex gradients of nonsmooth functions in direct search
  methods.
\newblock {\em IMA J. Numer. Anal.}, 28(4):770--784, 2008.

\bibitem{CustodioRochaVicente2010}
A.L. Cust{\'o}dio, H.~Rocha, and L.N. Vicente.
\newblock Incorporating minimum {F}robenius norm models in direct search.
\newblock {\em Comput. Optim. Appl.}, 46(2):265--278, 2010.

\bibitem{Gumma2014}
E.A.E. Gumma, M.H.A. Hashim, and M.M. Ali.
\newblock A derivative-free algorithm for linearly constrained optimization
  problems.
\newblock {\em Comput. Optim. Appl.}, 57(3):599--621, 2014.

\bibitem{W.L.Hare2014}
W.~Hare and Y.~Lucet.
\newblock Derivative-free optimization via proximal point methods.
\newblock {\em J. Optim. Theory App.}, 160(1):204--220, 2014.

\bibitem{hare-2010}
W.L. Hare.
\newblock Using derivative free optimization for constrained parameter
  selection in a home and community care forecasting model.
\newblock In {\em International Perspectives on Operations Research and Health
  Care, Proceedings of the 34th Meeting of the EURO Working Group on
  Operational Research Applied to Health Sciences}, pages 61--73, 2010.

\bibitem{CODE}
A.~Hedar.
\newblock Test functions for unconstrainted global optimization.
\newblock
  \url{http://www-optima.amp.i.kyoto-u.ac.jp/member/student/hedar/Hedar_files/TestGO_files/Page364.htm}.

\bibitem{Hiriart-Urruty-Lemarechal-1993a}
J.-B. Hiriart-Urruty and C.~Lemar{\'e}chal.
\newblock {\em Convex Analysis and Minimization Algorithms. {I}}, volume 305 of
  {\em Grundlehren der Mathematischen Wissenschaften [Fundamental Principles of
  Mathematical Sciences]}.
\newblock Springer-Verlag, Berlin, 1993.
\newblock Fundamentals.

\bibitem{Makela-2002}
M.M. M{\"a}kel{\"a}.
\newblock Survey of bundle methods for nonsmooth optimization.
\newblock {\em Optim. Methods Softw.}, 17(1):1--29, 2002.

\bibitem{marsden-feinstein-taylor-2008}
A.L. Marsden, J.A. Feinstein, and C.A. Taylor.
\newblock A computational framework for derivative-free optimization of
  cardiovascular geometries.
\newblock {\em Comput. Methods Appl. Mech. Engrg.}, 197(21-24):1890--1905,
  2008.

\bibitem{marsden-wang-dennis-moin-2007}
A.L. Marsden, M.~Wang, J.E. Dennis, Jr., and P.~Moin.
\newblock Trailing-edge noise reduction using derivative-free optimization and
  large-eddy simulation.
\newblock {\em J. Fluid Mech.}, 572:13--36, 2007.

\bibitem{More-Garbow-Hillstrom-1981}
J.J. Mor{\'e}, B.S. Garbow, and K.E. Hillstrom.
\newblock Testing unconstrained optimization software.
\newblock {\em ACM Trans. Math. Software}, 7(1):17--41, 1981.

\bibitem{MoreWild09}
J.J. Mor{\'e} and S.M. Wild.
\newblock Benchmarking derivative-free optimization algorithms.
\newblock {\em SIAM J. Optim.}, 20(1):172--191, 2009.

\bibitem{OliveiraSagastizabal2014}
W.~Oliveira and C.~Sagastiz\'abal.
\newblock Bundle methods in the {XXI}st century: A bird's-eye view.
\newblock {\em {Pesquisa Operacional}}, 34:647 -- 670, 12 2014.

\bibitem{Powell1992}
M.J.D. Powell.
\newblock A direct search optimization method that models the objective and
  constraint functions by linear interpolation.
\newblock In {\em Advances in optimization and numerical analysis ({O}axaca,
  1992)}, volume 275 of {\em Math. Appl.}, pages 51--67. Kluwer Acad. Publ.,
  Dordrecht, 1994.

\bibitem{Powell2002}
M.J.D. Powell.
\newblock Uobyqa: unconstrained optimization by quadratic approximation.
\newblock {\em Math. Program.}, 92(3):555--582, 2002.

\bibitem{Powell2004}
M.J.D. Powell.
\newblock Least {F}robenius norm updating of quadratic models that satisfy
  interpolation conditions.
\newblock {\em Math. Program.}, 100(1, Ser. B):183--215, 2004.

\bibitem{ScheinbergToint2010}
K.~Scheinberg and Ph.L. Toint.
\newblock Self-correcting geometry in model-based algorithms for
  derivative-free unconstrained optimization.
\newblock {\em SIAM J. Optim.}, 20(6):3512--3532, 2010.

\bibitem{ScottFrazierPowell2011}
W.~Scott, P.~Frazier, and W.~Powell.
\newblock The correlated knowledge gradient for simulation optimization of
  continuous parameters using {G}aussian process regression.
\newblock {\em SIAM J. Optim.}, 21(3):996--1026, 2011.

\bibitem{WildRegisShoemaker2008}
S.M. Wild, R.G. Regis, and C.A. Shoemaker.
\newblock {ORBIT}: optimization by radial basis function interpolation in
  trust-regions.
\newblock {\em SIAM J. Sci. Comput.}, 30(6):3197--3219, 2008.

\bibitem{YaoChenHuangvanTooren2014}
W.~Yao, X.Q. Chen, Y.Y. Huang, and M.~van Tooren.
\newblock A surrogate-based optimization method with {RBF} neural network
  enhanced by linear interpolation and hybrid infill strategy.
\newblock {\em Optim. Methods Softw.}, 29(2):406--429, 2014.

\end{thebibliography}
\bibliographystyle{plain} 

\newpage
\appendix

\section{Tables}

\begin{table}[ht] \caption{Results of each strategy of \DFPP~with no line search over all test problems.}\label{tab:nolineall}
\begin{center}
\begin{tabular}{|c|c|c|r|p{1.25cm}|c|c|c|r|}\cline{1-4}\cline{6-9}
$N_{\tt s}$	&	$N_{\tt n1}$	&	$N_{\tt n2}$	&	$\imp$	& &	$N_{\tt s}$	&	$N_{\tt n1}$	&	$N_{\tt n2}$	&	$\imp$	 \\\cline{1-4}\cline{6-9}
$\lin$	&	$\lin$	&	$\bilin$	&	249.01	&&	$\qua$	&	$\lin$	&	$\bilin$	&	118.98	\\
$\lin$	&	$\bilin$	&	$\lin$	&	247.74	&&	$\bilin$	&	$\mean$	&	$\qua$	&	118.09	\\
$\lin$	&	$\bilin$	&	$\bilin$	&	234.25	&&	$\mean$	&	$\mean$	&	$\lin$	&	117.51	\\
$\lin$	&	$\lin$	&	$\lin$	&	226.75	&&	$\mean$	&	$\mean$	&	$\bilin$	&	117.23	\\
$\lin$	&	$\mean$	&	$\lin$	&	212.73	&&	$\bilin$	&	$\bilin$	&	$\qua$	&	115.45	\\
$\lin$	&	$\mean$	&	$\bilin$	&	209.35	&&	$\bilin$	&	$\lin$	&	$\qua$	&	115.45	\\
$\lin$	&	$\qua$	&	$\lin$	&	205.58	&&	$\mean$	&	$\bilin$	&	$\bilin$	&	115.08	\\
$\lin$	&	$\lin$	&	$\mean$	&	204.31	&&	$\mean$	&	$\lin$	&	$\bilin$	&	115.08	\\
$\lin$	&	$\lin$	&	$\qua$	&	204.19	&&	$\qua$	&	$\mean$	&	$\bilin$	&	113.97	\\
$\lin$	&	$\mean$	&	$\mean$	&	202.25	&&	$\qua$	&	$\bilin$	&	$\mean$	&	113.80	\\
$\lin$	&	$\bilin$	&	$\qua$	&	202.06	&&	$\qua$	&	$\mean$	&	$\mean$	&	113.80	\\
$\bilin$	&	$\lin$	&	$\lin$	&	200.75	&&	$\qua$	&	$\lin$	&	$\mean$	&	113.80	\\
$\lin$	&	$\bilin$	&	$\mean$	&	199.24	&&	$\qua$	&	$\qua$	&	$\lin$	&	112.78	\\
$\lin$	&	$\mean$	&	$\qua$	&	198.53	&&	$\qua$	&	$\qua$	&	$\mean$	&	110.35	\\
$\lin$	&	$\qua$	&	$\mean$	&	197.41	&&	$\qua$	&	$\qua$	&	$\bilin$	&	110.23	\\
$\lin$	&	$\qua$	&	$\bilin$	&	195.55	&&	$\mean$	&	$\bilin$	&	$\mean$	&	108.27	\\
$\lin$	&	$\qua$	&	$\qua$	&	175.70	&&	$\mean$	&	$\mean$	&	$\mean$	&	108.27	\\
$\bilin$	&	$\bilin$	&	$\lin$	&	169.67	&&	$\mean$	&	$\lin$	&	$\mean$	&	108.27	\\
$\bilin$	&	$\mean$	&	$\lin$	&	143.23	&&	$\mean$	&	$\qua$	&	$\lin$	&	105.59	\\
$\bilin$	&	$\bilin$	&	$\bilin$	&	143.02	&&	$\bilin$	&	$\qua$	&	$\bilin$	&	103.94	\\
$\bilin$	&	$\lin$	&	$\bilin$	&	143.02	&&	$\bilin$	&	$\qua$	&	$\mean$	&	102.79	\\
$\mean$	&	$\lin$	&	$\lin$	&	137.59	&&	$\qua$	&	$\bilin$	&	$\qua$	&	102.45	\\
$\bilin$	&	$\mean$	&	$\bilin$	&	131.18	&&	$\qua$	&	$\qua$	&	$\qua$	&	102.45	\\
$\qua$	&	$\bilin$	&	$\lin$	&	130.88	&&	$\qua$	&	$\mean$	&	$\qua$	&	102.45	\\
$\qua$	&	$\mean$	&	$\lin$	&	130.83	&&	$\qua$	&	$\lin$	&	$\qua$	&	102.45	\\
$\bilin$	&	$\mean$	&	$\mean$	&	126.90	&&	$\mean$	&	$\bilin$	&	$\qua$	&	99.11	\\
$\bilin$	&	$\bilin$	&	$\mean$	&	126.10	&&	$\mean$	&	$\mean$	&	$\qua$	&	99.11	\\
$\bilin$	&	$\lin$	&	$\mean$	&	126.10	&&	$\mean$	&	$\lin$	&	$\qua$	&	99.11	\\
$\mean$	&	$\bilin$	&	$\lin$	&	124.79	&&	$\mean$	&	$\qua$	&	$\bilin$	&	94.82	\\
$\qua$	&	$\lin$	&	$\lin$	&	123.47	&&	$\mean$	&	$\qua$	&	$\mean$	&	91.17	\\
$\bilin$	&	$\qua$	&	$\lin$	&	121.25	&&	$\bilin$	&	$\qua$	&	$\qua$	&	91.08	\\
$\qua$	&	$\bilin$	&	$\bilin$	&	118.98	&&	$\mean$	&	$\qua$	&	$\qua$	&	84.84	\\
\cline{1-4}\cline{6-9}
\end{tabular}
\end{center}
\end{table}

\begin{table}[ht] \caption{Results of each strategy of \DFPP~with no line search for low dimensional problems.}\label{tab:nolinelow}
\begin{center}
\begin{tabular}{|c|c|c|r|p{1.25cm}|c|c|c|r|}\cline{1-4}\cline{6-9}
$N_{\tt s}$	&	$N_{\tt n1}$	&	$N_{\tt n2}$	&	$\imp$	& &	$N_{\tt s}$	&	$N_{\tt n1}$	&	$N_{\tt n2}$	&	$\imp$	 \\\cline{1-4}\cline{6-9}
$\lin$	&	$\qua$	&	$\lin$	&	148.88	   &&	$\qua$	&	$\bilin$	&	$\bilin$	&	78.55	\\
$\lin$	&	$\lin$	&	$\qua$	&	143.69	&&	$\qua$	&	$\lin$	&	$\mean$	&	77.74	\\
$\lin$	&	$\lin$	&	$\mean$	&	142.21	&&	$\qua$	&	$\bilin$	&	$\mean$	&	77.74	\\
$\lin$	&	$\lin$	&	$\bilin$	&	140.87	&&	$\qua$	&	$\mean$	&	$\mean$	&	77.74	\\
$\lin$	&	$\qua$	&	$\mean$	&	140.02	&&	$\mean$	&	$\qua$	&	$\lin$	&	77.65	\\
$\lin$	&	$\mean$	&	$\qua$	&	137.87	&&	$\bilin$	&	$\lin$	&	$\bilin$	&	76.91	\\
$\lin$	&	$\qua$	&	$\bilin$	&	135.96	&&	$\bilin$	&	$\bilin$	&	$\bilin$	&	76.91	\\
$\lin$	&	$\mean$	&	$\mean$	&	135.85	&&	$\bilin$	&	$\lin$	&	$\mean$	&	76.74	\\
$\lin$	&	$\lin$	&	$\lin$	&	134.79	&&	$\bilin$	&	$\bilin$	&	$\mean$	&	76.74	\\
$\lin$	&	$\bilin$&	$\lin$	    &	134.66	&&	$\mean$	&	$\lin$	&	$\mean$	&	76.05	\\
$\lin$	&	$\mean$	&	$\lin$	&	134.60	&&	$\mean$	&	$\bilin$	&	$\mean$	&	76.05	\\
$\lin$	&	$\mean$	&	$\bilin$	&	134.07	&&	$\mean$	&	$\mean$	&	$\mean$	&	76.05	\\
$\lin$	&	$\bilin$	&	$\qua$      &	131.24	&&	$\qua$	&	$\qua$	&	$\lin$	&	75.57	\\
$\lin$	&	$\bilin$	&	$\bilin$    &	131.15	&&	$\bilin$	&	$\mean$	&	$\mean$	&	75.24	\\
$\lin$	&	$\bilin$	&	$\mean$	&	127.92	&&	$\bilin$	&	$\mean$	&	$\qua$	&	74.75 \\
$\lin$	&	$\qua$	&	$\qua$	&	123.08	&&	$\qua$	&	$\qua$	&	$\bilin$	&	74.71	\\
$\bilin$	&	$\lin$	&	$\lin$	    &	109.17	&&	$\qua$	&	$\qua$	&	$\mean$	&	74.32	\\
$\mean$	&	$\lin$	&	$\lin$	&	102.88	&&	$\bilin$	&	$\lin$	&	$\qua$	&	71.73	\\
$\bilin$	&	$\bilin$	&	$\lin$      &	94.27	&&	$\bilin$	&	$\bilin$	&	$\qua$	&	71.73	\\
$\bilin$	&	$\mean$   	&	$\lin$	    &	93.35	&&	$\bilin$	&	$\qua$	&	$\bilin$	&	70.87	\\
$\qua$	&	$\mean$	&	$\lin$	&	93.10	&&	$\mean$	&	$\lin$	&	$\qua$	&	69.64 \\
$\qua$	&	$\bilin$	&	$\lin$	    &	93.09	&&	$\mean$	&	$\bilin$	&	$\qua$	&	69.64	\\
$\mean$	&	$\bilin$	&	$\lin$	    &	90.17	&&	$\mean$	&	$\mean$	&	$\qua$	&	69.64	\\
$\qua$	&	$\lin$	&	$\lin$	&	86.09	&&	$\bilin$	&	$\qua$	&	$\mean$	&	68.04	\\
$\bilin$	&	$\qua$	&	$\lin$	    &	85.88	&&	$\mean$	&	$\qua$	&	$\bilin$	&	66.80	\\
$\mean$	&	$\mean$	&	$\lin$	&	84.56	&&	$\qua$	&	$\lin$	&	$\qua$	&	64.68	\\
$\mean$	&	$\mean$	&	$\bilin$	&	80.79	&&	$\qua$	&	$\bilin$	&	$\qua$	&	64.68	\\
$\bilin$	&	$\mean$	&	$\bilin$	&	78.81	&&	$\qua$	&	$\mean$	&	$\qua$	&	64.68	\\
$\mean$	&	$\lin$	&	$\bilin$	&	78.71	&&	$\qua$	&	$\qua$	&	$\qua$	&	64.68	\\
$\mean$	&	$\bilin$	&	$\bilin$	&	78.71	&&	$\mean$	&	$\qua$	&	$\mean$	&	63.56	\\
$\qua$	&	$\mean$	&	$\bilin$	&	78.57	&&	$\bilin$	&	$\qua$	&	$\qua$	&	60.18	\\
$\qua$	&	$\lin$	&	$\bilin$	&	78.55	&&	$\mean$	&	$\qua$	&	$\qua$	&	55.04	\\
\cline{1-4}\cline{6-9}
\end{tabular}
\end{center}
\end{table}

\begin{table}[ht] \caption{Results of each strategy of \DFPP~with no line search for high dimensional problems.}\label{tab:nolinehigh}
\begin{center}
\begin{tabular}{|c|c|c|r|p{1.25cm}|c|c|c|r|}\cline{1-4}\cline{6-9}
$N_{\tt s}$	&	$N_{\tt n1}$	&	$N_{\tt n2}$	&	$\imp$	& &	$N_{\tt s}$	&	$N_{\tt n1}$	&	$N_{\tt n2}$	&	$\imp$	 \\\cline{1-4}\cline{6-9}
$\lin$	&	$\bilin$	&	$\lin$	&	113.09	&&	$\qua$	&	$\bilin$	&	$\qua$	&	37.77	\\
$\lin$	&	$\lin$	&	$\bilin$	&	108.14	&&	$\qua$	&	$\mean$	&	$\qua$	&	37.77	\\
$\lin$	&	$\bilin$	&	$\bilin$	&	103.09	&&	$\qua$	&	$\qua$	&	$\qua$	&	37.77	\\
$\lin$	&	$\lin$	&	$\lin$	&	91.96	&&	$\qua$	&	$\mean$	&	$\lin$	&	37.73	\\
$\bilin$	&	$\lin$	&	$\lin$	&	91.58	&&	$\qua$	&	$\lin$	&	$\lin$	&	37.38	\\
$\lin$	&	$\mean$	&	$\lin$	&	78.13	&&	$\qua$	&	$\qua$	&	$\lin$	&	37.21	\\
$\bilin$	&	$\bilin$	&	$\lin$	&	75.40	&&	$\mean$	&	$\mean$	&	$\bilin$	&	36.44	\\
$\lin$	&	$\mean$	&	$\bilin$	&	75.28	&&	$\mean$	&	$\lin$	&	$\bilin$	&	36.36	\\
$\lin$	&	$\bilin$	&	$\mean$	&	71.32	&&	$\mean$	&	$\bilin$	&	$\bilin$	&	36.36	\\
$\lin$	&	$\bilin$	&	$\qua$	&	70.82	&&	$\qua$	&	$\lin$	&	$\mean$	&	36.06	\\
$\lin$	&	$\mean$	&	$\mean$	&	66.39	&&	$\qua$	&	$\bilin$	&	$\mean$	&	36.06	\\
$\bilin$	&	$\lin$	&	$\bilin$	&	66.11	    &&	$\qua$	&	$\mean$	&	$\mean$	&	36.06	\\
$\bilin$	&	$\bilin$	&	$\bilin$	&	66.11	&&	$\qua$	&	$\qua$	&	$\mean$	&	36.02	\\
$\lin$	&	$\lin$	&	$\mean$	&	62.10	&&	$\qua$	&	$\qua$	&	$\bilin$	&	35.51	\\
$\lin$	&	$\mean$	&	$\qua$	&	60.66	&&	$\qua$	&	$\mean$	&	$\bilin$	&	35.40	\\
$\lin$	&	$\lin$	&	$\qua$	&	60.50	&&	$\bilin$	&	$\qua$	&	$\lin$	&	35.37	\\
$\lin$	&	$\qua$	&	$\bilin$	&	59.59	&&	$\bilin$	&	$\qua$	&	$\mean$	&	34.75	\\
$\lin$	&	$\qua$	&	$\mean$	&	57.39	&&	$\mean$	&	$\lin$	&	$\lin$	&	34.71	\\
$\lin$	&	$\qua$	&	$\lin$	&	56.70	&&	$\mean$	&	$\bilin$	&	$\lin$	&	34.62	\\
$\lin$	&	$\qua$	&	$\qua$	&	52.62	&&	$\bilin$	&	$\qua$	&	$\bilin$	&	33.07	\\
$\bilin$	&	$\mean$	&	$\bilin$	&	52.37	&&	$\mean$	&	$\mean$	&	$\lin$	&	32.96	\\
$\bilin$	&	$\mean$	&	$\mean$	&	51.66	&&	$\mean$	&	$\lin$	&	$\mean$	&	32.22	\\
$\bilin$	&	$\mean$	&	$\lin$	&	49.88	&&	$\mean$	&	$\bilin$	&	$\mean$	&	32.22	\\
$\bilin$	&	$\lin$	&	$\mean$	&	49.36	&&	$\mean$	&	$\mean$	&	$\mean$	&	32.22	\\
$\bilin$	&	$\bilin$	&	$\mean$	&	49.36	&&	$\bilin$	&	$\qua$	&	$\qua$	&	30.90	\\
$\bilin$	&	$\lin$	&	$\qua$	&	43.72	&&	$\mean$	&	$\qua$	&	$\qua$	&	29.80	\\
$\bilin$	&	$\bilin$	&	$\qua$	&	43.72	&&	$\mean$	&	$\lin$	&	$\qua$	&	29.47	\\
$\bilin$	&	$\mean$	&	$\qua$	&	43.33	&&	$\mean$	&	$\bilin$	&	$\qua$	&	29.47	\\
$\qua$	&	$\lin$	&	$\bilin$	&	40.43	&&	$\mean$	&	$\mean$	&	$\qua$	&	29.47	\\
$\qua$	&	$\bilin$	&	$\bilin$	&	40.43	&&	$\mean$	&	$\qua$	&	$\bilin$	&	28.03	\\
$\qua$	&	$\bilin$	&	$\lin$	&	37.78	&&	$\mean$	&	$\qua$	&	$\lin$	&	27.93	\\
$\qua$	&	$\lin$	&	$\qua$	&	37.77	&&	$\mean$	&	$\qua$	&	$\mean$	&	27.61	\\
\cline{1-4}\cline{6-9}
\end{tabular}
\end{center}
\end{table}

\begin{table}[ht] \caption{Results of each strategy of \DFPP~with line search over all test problems.}\label{tab:withlineall}
\begin{center}
\begin{tabular}{|c|c|c|r|p{1.25cm}|c|c|c|r|}\cline{1-4}\cline{6-9}
$N_{\tt s}$	&	$N_{\tt n1}$	&	$N_{\tt n2}$	&	$\imp$	& &	$N_{\tt s}$	&	$N_{\tt n1}$	&	$N_{\tt n2}$	&	$\imp$	 \\\cline{1-4}\cline{6-9}
$\lin$	&	$\lin$	&	$\lin$	&	377.91	&&	$\bilin$	&	$\qua$	&	$\lin$	&	237.10	\\
$\lin$	&	$\bilin$	&	$\lin$	&	305.56	&&	$\mean$	&	$\mean$	&	$\lin$	&	234.47	\\
$\lin$	&	$\lin$	&	$\bilin$	&	298.57	&&	$\qua$	&	$\bilin$	&	$\mean$	&	233.97	\\
$\lin$	&	$\qua$	&	$\lin$	&	297.85	&&	$\qua$	&	$\mean$	&	$\mean$	&	233.97	\\
$\lin$	&	$\mean$	&	$\lin$	&	296.60	&&	$\qua$	&	$\lin$	&	$\mean$	&	233.97	\\
$\lin$	&	$\lin$	&	$\qua$	&	282.40	&&	$\lin$	&	$\bilin$	&	$\mean$	&	233.34	\\
$\bilin$	&	$\lin$	&	$\lin$	&	282.27	&&	$\bilin$	&	$\mean$	&	$\qua$	&	232.09	\\
$\lin$	&	$\bilin$	&	$\qua$	&	279.63	&&	$\qua$	&	$\qua$	&	$\mean$	&	231.28	\\
$\lin$	&	$\qua$	&	$\bilin$	&	279.53	&&	$\mean$	&	$\bilin$	&	$\qua$	&	229.72	\\
$\lin$	&	$\bilin$	&	$\bilin$	&	278.18	&&	$\mean$	&	$\mean$	&	$\qua$	&	229.72	\\
$\lin$	&	$\mean$	&	$\qua$	&	271.84	&&	$\mean$	&	$\lin$	&	$\qua$	&	229.72	\\
$\bilin$	&	$\bilin$	&	$\lin$	&	267.72	&&	$\qua$	&	$\bilin$	&	$\qua$	&	228.39	\\
$\lin$	&	$\lin$	&	$\mean$	&	266.60	&&	$\qua$	&	$\qua$	&	$\qua$	&	228.39	\\
$\qua$	&	$\lin$	&	$\lin$	&	263.12	&&	$\qua$	&	$\mean$	&	$\qua$	&	228.39	\\
$\lin$	&	$\qua$	&	$\qua$	&	258.12	&&	$\qua$	&	$\lin$	&	$\qua$	&	228.39	\\
$\qua$	&	$\bilin$	&	$\lin$	&	256.75	&&	$\bilin$	&	$\qua$	&	$\qua$	&	227.36	\\
$\qua$	&	$\mean$	&	$\lin$	&	254.52	&&	$\mean$	&	$\qua$	&	$\lin$	&	226.71	\\
$\bilin$	&	$\mean$	&	$\lin$	&	251.04	&&	$\mean$	&	$\bilin$	&	$\bilin$	&	226.51	\\
$\lin$	&	$\mean$	&	$\bilin$	&	250.89	&&	$\mean$	&	$\lin$	&	$\bilin$	&	226.51	\\
$\qua$	&	$\qua$	&	$\lin$	&	246.53	&&	$\bilin$	&	$\mean$	&	$\bilin$	&	225.04	\\
$\bilin$	&	$\bilin$	&	$\bilin$	&	245.19	&&	$\mean$	&	$\mean$	&	$\bilin$	&	224.73	\\
$\bilin$	&	$\lin$	&	$\bilin$	&	245.19	&&	$\mean$	&	$\qua$	&	$\bilin$	&	223.34	\\
$\mean$	&	$\bilin$	&	$\lin$	&	242.61	&&	$\bilin$	&	$\bilin$	&	$\mean$	&	223.05	\\
$\mean$	&	$\lin$	&	$\lin$	&	242.11	&&	$\bilin$	&	$\lin$	&	$\mean$	&	223.05	\\
$\qua$	&	$\bilin$	&	$\bilin$	&	241.41	&&	$\bilin$	&	$\qua$	&	$\bilin$	&	222.63	\\
$\qua$	&	$\lin$	&	$\bilin$	&	241.41	&&	$\mean$	&	$\qua$	&	$\qua$	&	219.03	\\
$\qua$	&	$\mean$	&	$\bilin$	&	239.60	&&	$\bilin$	&	$\mean$	&	$\mean$	&	214.59	\\
$\bilin$	&	$\bilin$	&	$\qua$	&	238.97	&&	$\bilin$	&	$\qua$	&	$\mean$	&	214.35	\\
$\bilin$	&	$\lin$	&	$\qua$	&	238.97	&&	$\mean$	&	$\qua$	&	$\mean$	&	210.46	\\
$\qua$	&	$\qua$	&	$\bilin$	&	238.53	&&	$\mean$	&	$\bilin$	&	$\mean$	&	208.80	\\
$\lin$	&	$\qua$	&	$\mean$	&	237.97	&&	$\mean$	&	$\mean$	&	$\mean$	&	208.80	\\
$\lin$	&	$\mean$	&	$\mean$	&	237.64	&&	$\mean$	&	$\lin$	&	$\mean$	&	208.80	\\
\cline{1-4}\cline{6-9}
\end{tabular}
\end{center}
\end{table}

\begin{table}[ht] \caption{Results of each strategy of \DFPP~with line search for low dimensional problems.}\label{tab:withlinelow}
\begin{center}
\begin{tabular}{|c|c|c|r|p{1.25cm}|c|c|c|r|}\cline{1-4}\cline{6-9}
$N_{\tt s}$	&	$N_{\tt n1}$	&	$N_{\tt n2}$	&	$\imp$	& &	$N_{\tt s}$	&	$N_{\tt n1}$	&	$N_{\tt n2}$	&	$\imp$	 \\\cline{1-4}\cline{6-9}
$\lin$	&	$\lin$	&	$\lin$	&	153.38	&&	$\qua$	&	$\mean$	&	$\bilin$	&	128.62	\\
$\lin$	&	$\lin$	&	$\qua$	&	149.57	&&	$\qua$	&	$\qua$	&	$\lin$	&	128.58	\\
$\lin$	&	$\bilin$	&	$\qua$	&	147.67	&&	$\lin$	&	$\mean$	&	$\bilin$	&	128.43	\\
$\lin$	&	$\mean$	&	$\qua$	&	146.20	&&	$\qua$	&	$\lin$	&	$\qua$	&	127.88	\\
$\lin$	&	$\lin$	&	$\mean$	&	142.24	&&	$\qua$	&	$\bilin$	&	$\qua$	&	127.88	\\
$\lin$	&	$\qua$	&	$\lin$	&	141.89	&&	$\qua$	&	$\mean$	&	$\qua$	&	127.88	\\
$\qua$	&	$\lin$	&	$\lin$	&	141.23	&&	$\qua$	&	$\qua$	&	$\qua$	&	127.88	\\
$\bilin$	&	$\mean$	&	$\qua$	&	141.03	&&	$\qua$	&	$\qua$	&	$\mean$	&	127.86	\\
$\mean$	&	$\lin$	&	$\qua$	&	140.20	&&	$\lin$	&	$\bilin$	&	$\mean$	&	126.56	\\
$\mean$	&	$\bilin$	&	$\qua$	&	140.20	&&	$\lin$	&	$\mean$	&	$\mean$	&	125.88	\\
$\mean$	&	$\mean$	&	$\qua$	&	140.20   	&&	$\mean$	&	$\mean$	&	$\bilin$	&	125.31	\\
$\lin$	&	$\qua$	&	$\qua$	&	138.87	&&	$\mean$	&	$\lin$	&	$\bilin$	&	125.21	\\
$\bilin$	&	$\lin$	&	$\qua$	&	138.75	&&	$\mean$	&	$\bilin$	&	$\bilin$	&	125.21	\\
$\bilin$	&	$\bilin$	&	$\qua$	&	138.75	&&	$\qua$	&	$\qua$	&	$\bilin$	&	125.20	\\
$\lin$	&	$\lin$	&	$\bilin$	&	138.48	&&	$\lin$	&	$\bilin$	&	$\bilin$	&	124.94	\\
$\lin$	&	$\qua$	&	$\bilin$	&	138.35	&&	$\mean$	&	$\lin$	&	$\mean$	&	124.50	\\
$\lin$	&	$\qua$	&	$\mean$	&	137.88	&&	$\mean$	&	$\bilin$	&	$\mean$	&	124.50	\\
$\bilin$	&	$\lin$	&	$\lin$	&	136.62	&&	$\mean$	&	$\mean$	&	$\mean$	&	124.50	\\
$\mean$	&	$\lin$	&	$\lin$	&	136.49	&&	$\mean$	&	$\qua$	&	$\lin$	&	124.07	\\
$\qua$	&	$\mean$	&	$\lin$	&	136.46	&&	$\bilin$	&	$\qua$	&	$\lin$	&	123.75	\\
$\qua$	&	$\bilin$	&	$\lin$	&	136.32	&&	$\bilin$	&	$\lin$	&	$\mean$	&	123.05	\\
$\bilin$	&	$\bilin$	&	$\lin$	&	135.27	&&	$\bilin$	&	$\bilin$	&	$\mean$	&	123.05	\\
$\mean$	&	$\mean$	&	$\lin$	&	134.27	&&	$\bilin$	&	$\mean$	&	$\mean$	&	123.04	\\
$\bilin$	&	$\mean$	&	$\lin$	&	133.99	&&	$\bilin$	&	$\mean$	&	$\bilin$	&	121.54	\\
$\mean$	&	$\qua$	&	$\qua$	&	133.51	&&	$\lin$	&	$\bilin$	&	$\lin$	&	121.36	\\
$\mean$	&	$\bilin$	&	$\lin$	&	133.43  	&&	$\mean$	&	$\qua$	&	$\bilin$	&	120.69	\\
$\bilin$	&	$\qua$	&	$\qua$	&	130.36	&&	$\lin$	&	$\mean$	&	$\lin$	&	120.03	\\
$\qua$	&	$\lin$	&	$\bilin$	&	130.11	&&	$\bilin$	&	$\qua$	&	$\mean$	&	119.10	\\
$\qua$	&	$\bilin$	&	$\bilin$	&	130.11	&&	$\mean$	&	$\qua$	&	$\mean$	&	118.97	\\
$\qua$	&	$\lin$	&	$\mean$	&	130.05	&&	$\bilin$	&	$\lin$	&	$\bilin$	&	118.19	\\
$\qua$	&	$\bilin$	&	$\mean$	&	130.05	&&	$\bilin$	&	$\bilin$	&	$\bilin$	&	118.19	\\
$\qua$	&	$\mean$	&	$\mean$	&	130.05	&&	$\bilin$	&	$\qua$	&	$\bilin$	&	117.63	\\
\cline{1-4}\cline{6-9}
\end{tabular}
\end{center}
\end{table}

\begin{table}[ht] \caption{Results of each strategy of \DFPP~with line search for high dimensional problems.}\label{tab:withlinehigh}
\begin{center}
\begin{tabular}{|c|c|c|r|p{1.25cm}|c|c|c|r|}\cline{1-4}\cline{6-9}
$N_{\tt s}$	&	$N_{\tt n1}$	&	$N_{\tt n2}$	&	$\imp$	& &	$N_{\tt s}$	&	$N_{\tt n1}$	&	$N_{\tt n2}$	&	$\imp$	 \\\cline{1-4}\cline{6-9}
$\lin$	&	$\lin$	&	$\lin$	&	224.53	&&	$\qua$	&	$\lin$	&	$\mean$	&	103.92	\\
$\lin$	&	$\bilin$	&	$\lin$	&	184.19	&&	$\qua$	&	$\bilin$	&	$\mean$	&	103.92	\\
$\lin$	&	$\mean$	&	$\lin$	&	176.56	&&	$\qua$	&	$\mean$	&	$\mean$	&	103.92	\\
$\lin$	&	$\lin$	&	$\bilin$	&	160.09	&&	$\bilin$	&	$\mean$	&	$\bilin$	&	103.51	\\
$\lin$	&	$\qua$	&	$\lin$	&	155.96	&&	$\qua$	&	$\qua$	&	$\mean$	&	103.43	\\
$\lin$	&	$\bilin$	&	$\bilin$	&	153.24	&&	$\mean$	&	$\qua$	&	$\bilin$	&	102.65	\\
$\bilin$	&	$\lin$	&	$\lin$	&	145.65	&&	$\mean$	&	$\qua$	&	$\lin$	&	102.65	\\
$\lin$	&	$\qua$	&	$\bilin$	&	141.18	&&	$\mean$	&	$\lin$	&	$\bilin$	&	101.30	\\
$\lin$	&	$\lin$	&	$\qua$	&	132.83	&&	$\mean$	&	$\bilin$	&	$\bilin$	&	101.30	\\
$\bilin$	&	$\bilin$	&	$\lin$	&	132.45	&&	$\qua$	&	$\lin$	&	$\qua$	&	100.51  \\
$\lin$	&	$\bilin$	&	$\qua$	&	131.96	&&	$\qua$	&	$\bilin$	&	$\qua$	&	100.51	\\
$\bilin$	&	$\lin$	&	$\bilin$	&	127.00	&&	$\qua$	&	$\mean$	&	$\qua$	&	100.51	\\
$\bilin$	&	$\bilin$	&	$\bilin$	&	127.00   	&&	$\qua$	&	$\qua$	&	$\qua$	&	100.51	\\
$\lin$	&	$\mean$	&	$\qua$	&	125.63	&&	$\bilin$	&	$\lin$	&	$\qua$	&	100.23	\\
$\lin$	&	$\lin$	&	$\mean$	&	124.36	&&	$\bilin$	&	$\bilin$	&	$\qua$	&	100.23	\\
$\lin$	&	$\mean$	&	$\bilin$	&	122.46	&&	$\mean$	&	$\mean$	&	$\lin$	&	100.19	\\
$\qua$	&	$\lin$	&	$\lin$	&	121.89	&&	$\lin$	&	$\qua$	&	$\mean$	&	100.09	\\
$\qua$	&	$\bilin$	&	$\lin$	&	120.43	&&	$\bilin$	&	$\lin$	&	$\mean$	&	100.00	\\
$\lin$	&	$\qua$	&	$\qua$	&	119.25	&&	$\bilin$	&	$\bilin$	&	$\mean$	&	100.00	\\
$\qua$	&	$\mean$	&	$\lin$	&	118.06	&&	$\mean$	&	$\mean$	&	$\bilin$	&	99.42	\\
$\qua$	&	$\qua$	&	$\lin$	&	117.95	&&	$\bilin$	&	$\qua$	&	$\qua$	&	97.00	\\
$\bilin$	&	$\mean$	&	$\lin$	&	117.05	&&	$\bilin$	&	$\qua$	&	$\mean$	&	95.24	\\
$\bilin$	&	$\qua$	&	$\lin$	&	113.35	&&	$\bilin$	&	$\mean$	&	$\mean$	&	91.55	\\
$\qua$	&	$\qua$	&	$\bilin$	&	113.32  	&&	$\mean$	&	$\qua$	&	$\mean$	&	91.49	\\
$\lin$	&	$\mean$	&	$\mean$	&	111.75	&&	$\bilin$	&	$\mean$	&	$\qua$	&	91.06	\\
$\qua$	&	$\lin$	&	$\bilin$	&	111.30	&&	$\mean$	&	$\lin$	&	$\qua$	&	89.52	\\
$\qua$	&	$\bilin$    &	$\bilin$	&	111.30	&&	$\mean$	&	$\bilin$	&	$\qua$	&	89.52	\\
$\qua$	&	$\mean$	&	$\bilin$	&	110.99	&&	$\mean$	&	$\mean$	&	$\qua$	&	89.52  \\
$\mean$	&	$\bilin$	&	$\lin$	&	109.18	&&	$\mean$	&	$\qua$	&	$\qua$	&	85.52	\\
$\lin$	&	$\bilin$	&	$\mean$	&	106.78	&&	$\mean$	&	$\lin$	&	$\mean$	&	84.30	\\
$\mean$	&	$\lin$	&	$\lin$	&	105.62	&&	$\mean$	&	$\bilin$	&	$\mean$	&	84.30	\\
$\bilin$	&	$\qua$	&	$\bilin$	&	105.00	&&	$\mean$	&	$\mean$	&	$\mean$	&	84.30	\\
\cline{1-4}\cline{6-9}
\end{tabular}
\end{center}
\end{table}

\end{document}